\documentclass[12pt]{article}
\def\hind{\hangindent=2pc\hangafter=1}
\newfont{\smcaps}{cmcsc10 scaled\magstep1}

\newcommand{\E}{{\rm ~E\,}}
\newcommand{\MA}{{\rm ~MA\,}}
\newcommand{\ARMA}{{\rm ~ARMA\,}}
\newcommand{\ARIMA}{{\rm ~ARIMA\,}}

\newcommand{\AR}{{\rm ~AR\,}}

\newtheorem{thm}{THEOREM}

 at 12truept
\usepackage{epsfig}
\usepackage{latexsym}
\raggedright
\begin{document}

\title{PORTMANTEAU TESTS FOR ARMA MODELS WITH INFINITE VARIANCE}
\date{}
\author{{\smcaps By J.-W. Lin AND A.I. McLeod}\\
{\it The University of Western Ontario\/}
}

\maketitle
\bigskip
\hrule
\medskip
Jen-Wen Lin and A. Ian McLeod (2008).
Portmanteau Tests for ARMA Models with Infinite Variance.
{\it Journal of Time Series Analysis}, 29, 600-617

\newpage
\baselineskip=22pt
\hrule
\bigskip
{\bf Abstract.\/}

Autoregressive and moving-average (ARMA) models with stable Paretian
errors is one of the most studied models for time series with
infinite variance. Estimation methods for these models have been
studied by many researchers but the problem of diagnostic checking
fitted models has not been addressed. In this paper, we develop
portmanteau tests for checking randomness of a time series with
infinite variance and as a diagnostic tool for checking model
adequacy of fitted ARMA models. It is assumed that least-squares or
an asymptotically equivalent estimation method, such as Gaussian
maximum likelihood in the case of AR models, is used. And it is
assumed that the distribution of the innovations is IID stable
Paretian. It is seen via simulation that the proposed portmanteau
tests do not converge well to the corresponding limiting
distributions for practical series length so a Monte-Carlo test is
suggested. Simulation experiments show that the proposed test
procedure works effectively. Two illustrative applications to actual
data are provided to demonstrate that an incorrect conclusion may
result if the usual portmanteau test based on the finite variance
assumption is used.

\bigskip
{\bf Keywords.\/} ARMA models, Infinite variance, Least squares
method, Portmanteau test, Residual autocorrelation function, Stable
Paretian distribution

\hfill\eject

\begin{center}
1. INTRODUCTION
\end{center}
\medskip

Time series models with stable Paretian errors have been studied by
many researchers. Adler et al. (1998) discussed many aspects of how
to apply standard Box-Jenkins techniques to stable {\ARMA}
processes. Adler et al. (1998) concluded that, in principle, the
standard Box-Jenkins techniques do carry over to the stable setting
but a great deal of care needs to be exercised. In \S2 we briefly
review the stable Paretian distribution and in \S3 we develop
portmanteau tests for whiteness or randomness for an {\rm IID}
series. The whiteness test is illustrated with a brief application
to exchange rate data. In \S4 we develop portmanteau diagnostic
checks for residuals of an {\AR} model fitted by least-squares
assuming the true innovations are {\rm IID} stable Paretian
distributed. This is extended to the ARMA model in Appendix C. An
illustrative example shows the differences in inferences that may
result between the finite variance and infinite variance portmanteau
tests.

\begin{center}
2. THE STABLE PARETIAN DISTRIBUTION
\end{center}
\medskip

A stable distribution is usually defined through its characteristic function.
A random variable $Z$, or $Z_{\alpha }( \sigma ,\beta,\mu )$, is said to have a stable distribution
if its characteristic function has the following form:
\begin{equation}
{\E}(e^{i t Z })=\{
\begin{array}{c}
 \exp \left\{ -\sigma |t|^{\alpha}\, \left( 1-i \beta \ \ \mbox{sgn}(t)\, \tan \frac{\pi \alpha }{2}\right) +\,i \mu  t\right\} \ \ \ \ \ \mbox{if}
 \ \ \ \alpha \neq 1 \nonumber\\
 \exp \left\{ -\sigma  |t|\, \left( 1+i \beta  \frac{2}{\pi}\ \ \mbox{sgn}(t)\, \log |t| \right) +\,i \mu  t\right\}
\ \ \ \ \ \mbox{if} \ \ \ \alpha =1,
\end{array}
\newcounter{cf}
\setcounter{cf}{\value{equation}}
\end{equation}
where $i^2=-1$, $t$ is the parameter of the characteristic function,
$\alpha$ is the index of stability, or the characteristic
exponent, satisfying $0<\alpha \leq 2$, $\sigma >0$ is the scale
parameter, $\beta$ is the skewness satisfying $-1\leq \beta \leq
1$, $\mu\in R^1$ is the location parameter, and
\[
\mbox{sgn}(t)=\{
\begin{array}{c}
 1\ \ \ \mbox{ if}\ \ t >0 \\
 0\ \ \ \mbox{ if}\ \ t=0\\
 -1\ \ \mbox{if} \ \ t<0.
\end{array}
\]

In this paper, we restrict our attention to processes generated by
application of a linear filter to an independently and identically distributed
({\rm IID}) sequence, $\{Z_t: t=0,\pm 1,\ldots,\}$ , of random variables whose
distribution $F$ has Pareto-like tails, i.e.,
\begin{equation}
\{
\begin{array}{c}
x^\alpha\, (1-F(x))=x^\alpha\, P(Z_t >x)\rightarrow p\, C \\
x^\alpha\, F(-x)=x^\alpha\, P(Z_t <-x)\rightarrow q\, C,
\newcounter{Paret}
\setcounter{Paret}{\value{equation}}
\end{array}
\end{equation}
as $x\rightarrow \infty$, where
$0\leq\, p=1-q\leq\,1$, and $C$ is a finite positive constant,
or the dispersion of the random variable $Z_t$.

\bigskip
\begin{center}
3. PORTMANTEAU TESTS FOR RANDOMNESS OF STABLE PARETIAN TIME SERIES
\end{center}
\medskip

In this section, we shall derive the asymptotic distributions of
portmanteau tests for checking randomness of a sequence of stable
Paretian random variables. We consider the stable analogues of
portmanteau tests of Box and Pierce (1970) as well as Pe\v{n}a and
Rodriguez (2002), denoted by $Q_{BP}$ and $\hat D$, respectively. To
do so, we require some important properties of sample
autocorrelation functions (ACF) and sample partial autocorrelation
functions (PACF) of stable Paretian {\ARMA} processes (Brockwell and
Davis, 1991, Ch. 13; Samorodnitsky and Taqqu, 1994; Adler et al.,
1998).

\medskip
{\noindent \it 3.1 Asymptotic Distribution of Autocorrelation Function\hfill}
\medskip

Let $\{Z_t: t=0,\pm1,\pm2,\ldots\}$ be
an {\rm ~IID} sequence of
stable Paretian random variables
and $X_{t}$ be the strictly
stationary process defined by
\begin{equation}
X_{t}=\sum \limits_{j=-\infty }^{\infty }\psi _{j}Z_{t-j},\ \ t=1,\ldots
, n,
\end{equation}
where
\begin{equation}
\sum \limits_{j=-\infty }^{\infty } \left| j\right| \left|\psi _{j}\right|^{\delta } < \infty,
\ \ \mbox{for some} \ \ \delta \in \left( 0,\alpha \right) \cap \left[ 0,1\right] .
\end{equation}
The stable analogue of the autocorrelation function at lag $k$
is defined as
\begin{equation}
\rho_{k}=\sum \limits_{j}\psi _{j}\psi _{j+k} /\sum \limits_{j}\psi
_{j}^{2},\ \ k=1,2,\ldots.
\newcounter{defACF}
\setcounter{defACF}{\value{equation}}
\end{equation}
Eqn ({\thedefACF}) can be estimated by
the sample autocorrelation function as follows:
\begin{equation}
 {r}_{k}=\left\{ \sum \limits_{t =1}^{n-k} X_{t}X_{t+k}\right\}
/\sum \limits_{t =1}^{n} X_{t}^{2},\ \ k=1,2,\ldots ,
\newcounter{eACF}
\setcounter{eACF}{\value{equation}}
\end{equation}
for $\alpha\,>0$.
According to Davis and Resnick (1986), for any positive integer $k$,
the limiting distribution of sample autocorrelation functions is given by
\begin{equation}
\left[ \frac{n}{\log ( n) }\right] ^{\frac{1}{\alpha }}\left( r
_{1}-\rho_{1},\ldots ,{r}_{k}-\rho_{k}\right) ^{T}{\rightarrow
}\left( Y_{1},\ldots ,Y_{k}\right) ^{T},
\end{equation}
where  $\rightarrow$ denotes convergence in distribution and
\begin{equation}
Y_{h}=\sum \limits_{j=1}^{\infty }\,
\left( \rho_{k+j}+\rho_{k-j}-2 \rho_{j}\, \rho_{k}\right)
\,\frac{S_{j}}{S_{0}},\, h=1,\ldots ,k,
\end{equation}
where $S_{0},S_{1}, \ldots $ are independent stable variables; $S_{0}$
is positive with $S_{0}\sim Z_{\alpha /2}( C_{\alpha /2}^{-2/\alpha
},1,0) ,$ and the $S_{j}$ are $Z_{\alpha }( C_{\alpha }^{-1/\alpha
},0,0)$, where
\[
C_{\alpha }=
\frac{1-\alpha }{\Gamma ( 2-\alpha ) \cos ( \frac{\pi  \alpha }{2})
}\ \ \mbox{if} \ \ \alpha \neq 1,
\]
and
\[
C_{\alpha }=
\frac{2}{\pi }\ \ \mbox{if} \ \ \alpha =1.
\]

Under the null hypothesis that $X_{t}$ are a sequence of
{\rm~IID} stable Paretian random variables,
we have
$\rho_0=1$ and $\rho_k=0$ for $k \geq 1$
so
the limiting distribution of sample ACFs can be further simplified as follows:
\begin{equation}
\left[ \frac{n}{\log ( n) }\right] ^{\frac{1}{\alpha }}\,
\left( {r}_{1},\ldots ,{r}_{k} \right) ^{T}
{\rightarrow}
\,\left( W_{1},\ldots ,W_{k}\right) ^{T},
\newcounter{distSACF}
\setcounter{distSACF}{\value{equation}}
\end{equation}
where $W_{h}$ are given by
\begin{equation}
W_{h}=\frac{S_{h}}{S_{0}},\, h=1,\ldots , k.
\newcounter{Wdef}
\setcounter{Wdef}{\value{equation}}
\end{equation}
Note that, for $\alpha > 1$,
we may also use the mean-corrected sample autocorrelation function
 at lag $k$,
denoted as $\tilde r_k$, which is given by
\begin{equation}
 \tilde {r}_{k}=\sum \limits_{t =1}^{n-k}
 ( X_{t}-{\bar X})( X_{t+k}-{\bar X})
 /
\sum \limits_{t =1}^{n} ( X_{t}-{\bar X})^2,
\end{equation}
$k=1,2,\ldots.$
Davis and Resnick (1986) indicated that
the limiting distribution of $\tilde r_k$ is the same as that of $ r_k$.

\medskip
{\noindent \it 3.2 Asymptotic Distribution of Partial Autocorrelation Function\hfill}
\medskip

Consider an {\AR(p)} process,
\[
X_{t}-\phi _{1}X_{t-1}-\ldots -\phi _{p}X_{t-p}=Z_{t},
\]
where
$\{Z_t: t=0,\pm1,\pm2,\ldots\}$ are a sequence of {\rm~IID}
stable Paretian errors,
$1-\phi _{1}z-\ldots -\phi _{p}z^{p}\,\neq 0$, $|z|\,\leq 1$.
 Let $\rho_{(p)}=(\rho_{1},\ldots ,\rho_{p})^{T}$ be a vector of
 autocorrelation functions,
$\mathcal{R}_{(p)}=(\rho_{|i-j|})_{p\times p}$ be the $p\times p$
autocorrelation matrix, and $\phi _{(p)}=(\phi _{1},\ldots ,\phi
_{p})^{T}$.
The Yule-Walker equations are defined as
\begin{equation}
\mathcal{R}_{\left( p\right) }\phi _{\left( p\right) }=\rho_{\left( p\right)
}.
\end{equation}
The PACF at lag $p$ is simply
the $p$-th element of the solution of the Yule-walker equations,
\[
\phi _{(p)}^{YW}=\Psi\left(\rho_{(p)}\right)=\mathcal{R}_{(p)}^{-1}\rho_{(p)}.
\]
Likewise, the sample partial autocorrelation function at lag $p$ is defined as the $p$-th element of
the sample estimate of the Yule-walker solution,
\[
\hat{\phi }_{(p)}^{YW}=\Psi(\textbf{r}_{(p)})={\textbf{R}}_{(p)}^{-1} {\textbf{r}}_{(p)},
\]
where ${\textbf{R}}_{(p)}=(r_{|i-j|})_{p\times p}$
and ${\textbf{r}}_{(p)}=(r_1,\ldots,r_p)^T$ are
the $p\times p$ sample autocorrelation matrix and
the $p\times 1$ vector of sample autocorrelation functions, respectively.
It is apparent that
the sample partial autocorrelations is a function of sample autocorrelations.
Their relationship is clearly described in the Durbin-Levison algorithm.

Let ${\pi }_{k}$ be the sample PACF at lag $k$, and
${\pi }_{(m)}=({\pi}_{1},\ldots ,{\pi }_{m})^{T}$.
By the Durbin-Levison algorithm,
the vector ${\pi }_{(m)}$ can be expressed as a function
of $\textbf{r}_{(m)}$, ${\pi }_{(m)}=\psi ( \textbf{r}_{(m)})
$, with the $k$-th element given by
\begin{equation}
{\pi }_{k}=\psi ( \textbf{r}_{\left( k\right) })
=\frac{r_{k}-\textbf{r}_{\left( k-1\right) }^{T} \textbf{R}_{\left( k-1\right) }^{-1}
\textbf{r}_{\left( k-1\right) }^{*}}
{
1-\textbf{r}_{\left( k-1\right) }^{T} \textbf{R}_{\left(
k-1\right) }^{-1} \textbf{r}_{\left( k-1\right) }
},
\newcounter{DurbL}
\setcounter{DurbL}{\value{equation}}
\end{equation}
where $\textbf{R}_{(k)}$ and $\textbf{r}_{(k)}$ are as defined above
and $\textbf{r}_{(k)}^{*}=(r_{k},\ldots ,r_{1})^{T}$.

Following the proof in Monti (1994), we can derive
the asymptotic distribution of sample partial autocorrelation functions.
Under the null hypothesis that $X_{t}$ are independent,
the autocorrelation functions are all zero, and
according to Brockwell and Davis (1991, ch. 13),
\[
r_h =O_{p}\left( \left[ \frac{n}{\log ( n) }\right]^{-1/\alpha}\right),\ \ h=1,2, \ldots.
\]
Therefore,
\[
\textbf{R}_{(k)}=\textbf{1}_{k}+O_{p}\left( \left[\frac{n}{\log ( n) }\right]^{-1/\alpha}\right),
\]
where $\textbf{1}_k$ is a $k\times\, k$ identity matrix.
By eqn. ({\theDurbL}),
\begin{equation}
{\pi }_{\left( m\right) }
=
\textbf{r}_{\left( m\right) }+O_{p}\left(
\left[ \frac{n}{\log ( n) }\right] ^{-2/\alpha }\right).
\newcounter{phiandrhat}
\setcounter{phiandrhat}{\value{equation}}
\end{equation}
Using eqn. ({\thedistSACF}), we have
\begin{equation}
 \left[ \frac{n}{\log ( n) }\right] ^{\frac{1}{\alpha }}\,\left( \pi_{1},\ldots ,\pi _{m}\right) ^{T}
 {\rightarrow }\,
 \left(W_{1},\ldots ,W_{m}\right) ^{T}.
\newcounter{distPACF}
\setcounter{distPACF}{\value{equation}}
\end{equation}

\medskip
{\noindent \it 3.3 Asymptotic Distributions of $Q_{\rm BP}$ and $\hat D$ Tests\hfill}
\medskip

We can now derive the limiting distributions
of the $Q_{\rm BP}$ and $\hat D$ tests for checking randomness
of a sequence of stable Paretian random variables.
Under the assumption that $1<\alpha < 2$,
 Runde (1997)
derived the limiting distribution of $Q_{\rm BP}$,
based on the mean corrected sample autocorrelation functions.
His result is given by
\begin{equation}
\left( \frac{n}{\log ( n) }\right) ^{2/\alpha }\sum \limits_{j=1}^{m}
\tilde {r}_{j}^{2}\rightarrow W_{1}^{2}+\cdots +W_{m}^{2},
\newcounter{QBP}
\setcounter{QBP}{\value{equation}}
\end{equation}
where $\{W_k : k=1,\ldots,m\}$ are defined in eqn. ({\theWdef}).
Note that if $0< \alpha \leq1$, the limiting distribution of eqn. ({\theQBP})
remains the same if $\tilde r_k$ are replaced by $ r_k$.

Consider next the $\hat D$ test of Pe\v{n}a and Rodriguez (2002).
The test statistic may be given by
\begin{equation}
\hat{D}=\left(\frac{n}{\log(n)}\right)^{2/\alpha}\left ( 1-|\textbf{R}_{\left( m\right) }| ^{1/m}\right).
\newcounter{ddhat}
\setcounter{ddhat}{\value{equation}}
\end{equation}
Following the proof of Theorem 1 in Pe\v{n}a and Rodriguez (2002),
we may have the asymptotic distribution of eqn. ({\theddhat}) in the
following Theorem. The proof is given in Appendix A.

\begin{thm}
$\hat{D}$ in eqn. ({\theddhat}) is asymptotically distributed as
\[
\sum^m_{i=1} \frac{m+1-i}{m}\, W^2_{i},
\]
where $\{W_i: i=1,\ldots, m\}$
are as defined in eqn. ({\theWdef}).
\end{thm}

\medskip
{\noindent \bf Remark 1:}
It is possible to compute the
limiting distributions of
the $Q_{\rm BP}$ and $\hat D$ tests
by making use of the change variable technique and
some numerical algorithms of
calculating the probability density function
of stable random variables,
such as Mittnik et al. (1999).
This approach requires, however, intensive numerical computations.

\medskip
{\noindent \bf Remark 2:}
Another approach to obtaining the asymptotic distributions
of the $Q_{\rm BP}$ and $\hat D$ tests
is to simulate the aforementioned tests based on
their asymptotic distributions. For example,
$\hat D$ is simulated as defined in Theorem 1.
This approach also requires a large scale of computation but
is much less intensive computationally than the approach
mentioned in Remark 1.
This approach will be adopted in the subsequent analysis
based on
$10^4$ simulations.

\medskip
{\noindent \it 3.4 Simulation Experiments \hfill}
\medskip

The finite sample performance of $Q_{\rm BP}$ and $\hat D$ tests for
randomness will be investigated in this section. Based on $250$
simulations, the 5, 10, 30, 50, 70, 90, 95, 97.5, 99 $(\%)$
empirical quantiles of both tests with lag $m=5$ were calculated and
plotted against the corresponding asymptotic distributions. It is
seen in Figure 1 and Figure 2 that the empirical and asymptotic
quantiles do not agree very well unless $n$ is very large.

It is seen in Figures 1 to 2 that the speed of convergence
of both tests to the corresponding asymptotic distributions
 is very slow.
A solution to this problem is to use the Monte-Carlo test or
parametric bootstrap (Appendix B).

\begin{center}
[Figures 1 and 2 about here]
\end{center}

Consider the simulation experiments. {\rm IID} random sequence of
$Z_\alpha(1,0,0)$ with series length $n=250$ and $\alpha=$ $1.9$,
$1.7$, $1.5$, $1.3$, $1.1$ were simulated. The empirical sizes of
both tests were calculated based on $N=10^4$ simulations and each
Monte-Carlo test was simulated based on $10^3$ simulations. The
results are tabulated in Table 1. It is seen that the empirical
sizes of both tests are very close to the 5\% nominal level even
with $n=250$.

\begin{center}
[Table 1 about here]
\end{center}

\medskip
{\noindent \it 3.5 Illustrative Example \hfill}
\medskip

Consider the daily {\rm Canada/U.S.} exchange rates dated from
September 06, 1996 to September 05, 2006. The data was retrieved
from the website  of the Federal Reserve Bank of St. Louis and the
returns, $e_t = \log(z_{t+1}/z_t)$, were computed and tested for
randomness. The consistent estimators of McCulloch (1986) were used
to estimate $\alpha$ and $\beta$ for the returns. We obtained $\hat
\alpha_M=1.5644$ and $\hat \beta_M=-0.0472$. It is seen that $\hat
\beta_M$ is close to zero so the series is not highly skewed. Since
$\hat \alpha_M$ is much less than $2$, the usage of the portmanteau
tests in \S3 are more reasonable than that of the ordinary
portmanteau tests in this data. The P-values for $Q_{\rm LB}(m)$
test were determined using the asymptotic $\chi^2(m)$ distribution
and the Monte-Carlo method in Appendix B. The results are compared
in Table 2. Note that when $m=5$ the finite-variance portmanteau
test suggested possible evidence of non-randomness but this is not
the case when the infinite-variance Monte Carlo test is used.

\begin{center}
[Table 2 about here]
\end{center}

\medskip
{\noindent \bf Remark 3:} Portmanteau tests based on the
nonparametric bootstrap procedure could also be used but it would be
expected that they would be less powerful since less information is
used.

\bigskip
\begin{center}
4. DIAGNOSTIC CHECK FOR MODEL ADEQUACY OF AR$(p)$ MODELS WITH STABLE PARETIAN ERRORS
\end{center}
\medskip

\medskip
{\noindent \it 4.1 Some Asymptotic Results \hfill}
\medskip

In this section, we shall derive the asymptotic
distributions of $Q_{\rm BP}$ and $\hat D$ tests
for diagnostic check in model adequacy of {\AR}($p$) models
with stable Paretian errors.
Consider the general {\AR}($p$) process as follows:
\begin{equation}
 \phi(B)X_t=Z_t,
\newcounter{ARp}
\setcounter{ARp}{\value{equation}}
\end{equation}
where $\{Z_t: t=0,\pm 1, \pm 2,\ldots\}$ is an {\rm~IID} sequence of
stable Paretian random variables, $B$ denotes the backward operator,
and $\phi(B)=1-\phi_1 B-\cdots-\phi_p B^p$.
Let $\hat \phi_{(p)}=(\hat \phi_1,\ldots,\hat \phi_p)^T$ denote
the estimates of autoregressive coefficients.
The residuals of the fitted model are given as follows:
\begin{equation}
\hat Z_t=Z_t(\hat \phi_{(p)})=X_t-\hat \phi_1 X_{t-1}-\ldots-\hat \phi_p X_{t-p}=\hat \phi(B) X_t,
\newcounter{resid}
\setcounter{resid}{\value{equation}}
\end{equation}
and the corresponding residual autocorrelation at lag $k$ is given by
\[
\hat r_k=\frac{\sum \hat Z_t \hat Z_{t-k}}{\sum \hat Z^2_t}.
\]

Consider the estimators of $\hat \phi_{(p)}$ satisfying
\[
\hat \phi_{(p)}=\phi_{(p)}+O_p\left([n/\log(n)]^{-1/\alpha}\right).
\]
From Appendix C, the residual autocorrelation at lag $k$, $\hat
r_k$, can be approximated by the first order Taylor expansion about
error autocorrelation functions, $r_k$. Specifically, the
approximation is
\begin{equation}
\hat r_k=r_k+\sum_{j=1}^p (\phi_j-\hat \phi_j)\,\psi_{k-j}
+O_p\left([n/\log(n)]^{-2/\alpha}\right),
\newcounter{chthreeexpansion}
\setcounter{chthreeexpansion}{\value{equation}}
\end{equation}
where $\psi_j$ is the impulse response coefficient at lag $j$ and
$r_k=\sum Z_t Z_{t-k}/\sum Z^2_t$ is the error autocorrelation at lag $k$.
Eqn. ({\thechthreeexpansion}) can also be
written in matrix form, to order $O_p\left([n/\log(n)]^{-2/\alpha}\right)$,
\begin{equation}
\hat \textbf{r}_{(p)}= \textbf{r}_{(p)}+\textbf{X}\left( \phi_{(p)}-\hat \phi_{(p)}\right),
\newcounter{matrixform}
\setcounter{matrixform}{\value{equation}}
\end{equation}
where
\begin{equation}
\textbf{X}=\left[
\begin{array} {cccc}
1      & 0       & \cdots & 0 \\
\psi_1 & 1       & \ddots & 0 \\
\vdots & \vdots  & \ddots & 0 \\
\vdots & \vdots  & \ddots & 0 \\
\psi_{m-1} & \psi_{m-2} & \cdots & \psi_{m-p}
\end{array}
\right].
\newcounter{matrixcoef}
\setcounter{matrixcoef}{\value{equation}}
\end{equation}

By making use of eqn. ({\thechthreeexpansion}) or eqn. ({\thematrixform})
as well as following the proof in Theorem 1, we may derive
the asymptotic distributions of the aforementioned portmanteau tests
for diagnostic check in {\AR}$(p)$ models.
This distribution, however, is usually very complicated
and may not be traceable unless the {\AR}$(p)$
models of interest are fitted by least squares ({\rm LS}).
For simplicity, we only consider the case that eqn. ({\theARp}) is
estimated using least squares in the subsequent analysis.

According to $\S$4 in Davis (1996), if the
{\ARMA} parameters, $\beta$, are estimated using
least squares , we have
$ [n/{\log(n)}]^{1/\alpha}\left(\hat \beta_{LS}-\beta\right) $
converges in distribution, where $\hat \beta_{LS}$ denotes
the {\rm LS} estimates of $\beta$. Hence, in terms of our notation,
we have
$\hat \phi_{(p)}-\phi_{(p)}
=O_p\left([n/{\log(n)}]^{-1/\alpha}\right)$.
Then, by Box and Pierce (1970),
$\{\hat Z_t\}$ in eqn. {(\theresid)} satisfy the orthogonality conditions
and, to order $O_p\left(1/\sqrt{n}\, \, [n/{\log(n)}]^{-1/\alpha}\right)$,
\begin{equation}
\hat \textbf{r}_{(p)}^T\,\textbf{X}=0.
\newcounter{oc}
\setcounter{oc}{\value{equation}}
\end{equation}
If we now multiply eqn. ({\thematrixform}) on both sizes by
\[
\textbf{Q}=\textbf{X}(\textbf{X}^T\textbf{ X})^{-1}\textbf{X}^T,
\]
then using eqn. {(\theoc)} we have
\begin{equation}
\hat \textbf{r}_{(p)}=(\textbf{1}_m-\textbf{Q})\, \textbf{r}_{(p)}
\newcounter{simracf}
\setcounter{simracf}{\value{equation}}
\end{equation}
approximately, where $\textbf{1}_m$ is an $m \times m$ identity matrix and
$\textbf{Q}=\textbf{X}(\textbf{X}^T \textbf{X})^{-1}\textbf{X}^T$.
It was shown by Box and Pierce (1970) that
$\textbf{1}_m-\textbf{Q}$ is idempotent
of rank $m-p$.
Hence, the asymptotic distribution of the $Q_{\rm BP}$ test
is given by
\begin{equation}
(\frac{n}{\log n})^{2/\alpha} \sum_1^m \hat{r}_k^2 \rightarrow
\textbf{W}_m^T(\textbf{1}_m-\textbf{Q})\textbf{W}_m,
\newcounter{distQBP}
\setcounter{distQBP}{\value{equation}}
\end{equation}
where $\textbf{W}_m=(W_1,\ldots,W_m)^T$ and $\{W_i: i=1,\ldots,m\}$
are defined in eqn. ({\theWdef}).

Consider next the asymptotic distributions of residual partial autocorrelations.
Let $\hat \pi_{(m)}$ be the vector of the first $m$ residual partial autocorrelations
and $\pi_{(m)}$ is the vector of error partial autocorrelations.
The Taylor expansion of $\psi (\hat \textbf{r}_{(m)})$
around $ \textbf{r}_{(m)}$ yields
\begin{equation}
\hat \pi_{(m)}= \pi_{(m)}+\frac{\partial \pi_{(m)}}{\partial \textbf{r}_{(m)}}
\left(\hat \textbf{r}_{(m)}- \textbf{r}_{(m)}\right)
+O_p\left(\left[\frac{n}{\log n}\right]^{-2/\alpha}\right).
\newcounter{expn}
\setcounter{expn}{\value{equation}}
\end{equation}
By eqn. ({\theDurbL}) and ({\thephiandrhat}),  eqn. ({\theexpn}) becomes
\begin{equation}
\hat \pi_{(m)}=\hat \textbf{r}_{(m)}+O_p\left(\left[\frac{n}{\log n}\right]^{-2/\alpha}\right).
\newcounter{distpidd}
\setcounter{distpidd}{\value{equation}}
\end{equation}

Consider the Pe\v{n}a-Rodriguez test as the form of
\begin{equation}
\hat D =(\frac{n}{\log n})^{2/\alpha} \left(1-|\hat \textbf{R}_{(m)}|^{1/m}\right),
\newcounter{PRtest}
\setcounter{PRtest}{\value{equation}}
\end{equation}
where $\hat \textbf{R}_{(m)}=(\hat r_{|i-j|})_{m,m}$
is the $m\times m$ residual autocorrelation matrix.
By eqn. ({\thedistpidd}) and following
the proof in Theorem 1,
the limiting distribution of
eqn. ({\thePRtest}) is $\textbf{W}_m^T\, \textbf{A}_m\,\textbf{W}_m$,
where
$\textbf{A}_m=(\textbf{1}_m-\textbf{Q})^T\,\mathcal{W}_{m,m}\,
(\textbf{1}_m-\textbf{Q})$ and
$\mathcal{W}_{m,m}$ is a $m \times m$ diagonal matrix with $(i,i)$-th element equal to $(m-i+1)/m$ for
$i=1,\cdots,m$.

\medskip
{\noindent \bf Remark 4:} It is shown in Appendix C.4 that the
residuals in a fitted ARMA model are asymptotically equivalent to
those in a particular AR model. Hence the asympotic results for the
AR may be extended to the ARMA case.

\medskip
{\noindent \it 4.2 Some Size and Power Calculations \hfill}
\medskip

As in $\S$3.4, the slow convergence of $Q_{\rm BP}$ and $\hat D$ tests
to their asymptotic distributions is also present at the residual
autocorrelations. The first order autoregressive process $X_t=0.5
X_{t-1}+Z_t$ with $Z_t \sim Z_{1.2}(1,0,0)$ was simulated and
{\AR}$(1)$ models were fitted to the data. Then the 5, 10, 30, 50,
70, 90, 95, 97.5, 99 $(\%)$ empirical quantiles of $\hat r_1$ were
plotted against its theoretical asymptotic distribution based on
$10^3$ simulations. The asymptotic distribution of the error
autocorrelation at lag one, $ r_1$, was also plotted in Figure 3. It
is seen that empirical quantiles of $\hat r_1$ get closer to its
asymptotic distribution as the series length $n$ increases. However,
this is not the case for the empirical quantiles of $\hat r_1$ to
the asymptotic distribution of $r_1$. Therefore, serious size
distortion may be present in this case if one uses error
autocorrelations as a diagnostic tool for checking model adequacy.
The slow convergence of residual autocorrelations to its asymptotic
distribution may cause difficulties in using portmanteau tests in
practice. Therefore, as in $\S$3.4, we suggested using the
Monte-Carlo test to improve the effectiveness of portmanteau tests.

\begin{center}
[Figure 3]
\end{center}

We now investigate the effectiveness of $Q_{\rm BP}$ and $\hat D$ tests
for diagnostic check in fitted {\AR} models with stable Paretian
errors. The empirical sizes of $\hat D$ and $Q_{\rm BP}$ tests for a
$5\%$ significance test were first calculated via simulation. In
this experiment, {\AR}$(1)$ models, $X_t=\phi_1 X_{t-1}+Z_t$, were
simulated, where $Z_t\sim Z_{1.5}(1,0,0)$ and $\phi_1=0,\pm 0.1,\pm
0.3,\pm 0.5,\pm 0.7,\pm 0.9$ and {\AR}$(1)$ models were fitted to
the simulated data by the Burg algorithm. The empirical size for
each test was calculated based on $N=10^4$ simulations and each
Monte Carlo test used $10^3$ simulations. Series length $n=100$ and
lags $m=5,10,20$ were investigated. It is seen in Table 3 that the
empirical sizes of both tests are very close to their nominal level.

\begin{center}
[Table 3]
\end{center}

The empirical powers of $\hat D$ and $Q_{\rm BP}$ tests as
diagnostic tools were also investigated via simulation. Twelve
{\ARMA}$(2,2)$ models of series length $n=100$ in Table 4 of
Pe\v{n}a and Rodriguez (2002) were simulated and {\AR}$(1)$ models
were fitted to the simulated data using the Burg algorithm. Both
tests with lags $m=5,10,20$ were calculated using the parametric
bootstrap procedure. The empirical powers were calculated based on
$N=10^3$ simulations and each Monte Carlo test used $10^3$
simulations. It is seen in Table 4 that the empirical powers of both
tests are reasonably good for most models. Some of them are even
better than the powers listed in Pe\v{n}a and Rodriguez (2002). In
addition, increasing the series length can also improve the
effectiveness of the proposed test procedure. For example, with
model 3 in Table 2, if the series length was increased to $n=250$,
the empirical powers of the $\hat D$ test at lags $m=5,10,20$ were
increased significantly from 23.37\%, 20.10\% and 17.61\% to
58.27\%, 43.71\% and 35.52\%, respectively. Similar improvement was
also found in the $Q_{\rm BP}$ test. Finally, as in Pe\v{n}a and
Rodriguez (2002), our simulation experiments show that $\hat D$ is
more powerful than $Q_{\rm BP}$ as a diagnostic tool.

\begin{center}
[Table 4]
\end{center}

\medskip
{\noindent \bf Remark 5:} It is well known that the Burg estimate of
$\phi_1$ is close to the {\rm LS} estimate. The advantage of using
Burg estimate is that it is always in the stationary region and this is needed for
the Monte-Carlo test.

\medskip
{\noindent \it 4.3 Illustrative Application \hfill}
\medskip

Tsay (2002, Ch. 2) tentatively identified an AR(3) or AR(5) model
for the monthly simple returns of CRSP value-weighted index from
January 1926 to December 1997 using the partial autocorrelation
function. Here $n=864$ and the usual Box-Pierce portmanteau test at
lags $m=5,10,20$ does not suggest model inadequacy of either model
at the 5\% level. By applying our Monte-Carlo test procedure,
however, both the $\hat D$ and $Q_{\rm BP}$ tests in \S4 reject both
models.
The P-values are displayed in Table 5.
 The infinite variance hypothesis is plausible since the estimates
for $\alpha$ of residuals in the fitted AR$(3)$ and AR$(5)$ models
are 1.696 and 1.635, respectively. We may conclude from this example
that using the ordinary portmanteau tests may lead to a wrong
decision if innovations have infinite variance.

\begin{center}
[Table 5]
\end{center}

\bigskip
\begin{center}
5. CONCLUDING REMARK
\end{center}
\medskip

We will provide an R package implementing the portmanteau tests
described in this paper on CRAN.

\newpage
\bigskip
\begin{center}
APPENDIX A: PROOF OF THEOREM 1
\end{center}
\medskip

First, by decomposing the determinant of the sample autocorrelation matrix $\textbf{R}_{(m)}$,
Pena and Rodriguez (2002) showed that
$|\textbf{R}_{(m)}|^{1/m}$ is a weighted function of the first $m$ partial autocorrelations.
Specifically,
\begin{equation}
|\textbf{R}_{(m)}|^{1/m}=\prod^m_{i=1} (1-\pi^2_i)^{(m+1-i)/m}.
\newcounter{detdec}
\setcounter{detdec}{\value{equation}}
\end{equation}
Suppose that under the null hypothesis, $\hat D$
is asymptotic distributed as $\mathcal{X}$.
By applying the $\delta$-method to $g(x)=\log(1-x)$, it follows that
$-\left(n / \log(n)\right)^{2/\alpha}\, \log \left(|\textbf{R}_{(m)}|^{1/m}\right)$
is asymptotically distributed as $\mathcal{X}$.
From eqn. ({\thedetdec}), we can have 
\begin{eqnarray}
&-&\left(\frac{n}{\log(n)}\right)^{2/\alpha} \log \left(|\textbf{R}_{m}|^{1/m}\right)=
\nonumber \\
&-&\left(\frac{n}{\log(n)}\right)^{2/\alpha} \sum \limits_{i=1}^{m}\frac{m-i+1}{m}\log( 1-{\pi }_{i}^{2}).
\end{eqnarray}
Next suppose that
\begin{equation}
\left(\frac{n}{\log(n)}\right)^{2/\alpha}\,
\left({\pi }_{1}^{2},
{\pi }_{2}^{2},
\ldots,
{\pi }_{m}^{2}\right)^{T}\longrightarrow Y,
\newcounter{AT}
\setcounter{AT}{\value{equation}}
\end{equation}
and apply the multivariate $\delta$-method to
\[
g( {\pi }_{1}^{2},{\pi }_{2}^{2},\ldots ,{\pi }_{m}^{2})
=-\sum \limits_{i=1}^{m}\frac{m-i+1}{m}\log (1-{\pi }_{i}^{2}),
\]
it follows that
\begin{equation}
-\sum \limits_{i=1}^{m}\frac{m-i+1}{m}\log (1-{\pi }_{i}^{2})
\rightarrow
\left( 1,\frac{m-1}{m},\ldots,\frac{1}{m}\right)Y.
\end{equation}
From the Cramer-Wold theorem, it follows that
\begin{eqnarray}
\left( 1,\frac{m-1}{m},\cdots,\frac{1}{m}\right)
\left( \left(\frac{n}{\log(n)}\right)^{2/\alpha} \pi_{1}^{2}
,\ldots , \left(\frac{n}{\log(n)}\right)^{2/\alpha}\pi_{m}^{2}\right)
^{T} \nonumber \\
\longrightarrow \left( 1,\frac{m-1}{m},\ldots ,\frac{1}{m}\right)
Y
\newcounter{ATT}
\setcounter{ATT}{\value{equation}}
\end{eqnarray}
By eqn. ({\thedistPACF}), it follows that
\begin{eqnarray}
\left( 1,\frac{m-1}{m},\ldots,\frac{1}{m}\right)
\left( \left(\frac{n}{\log(n)}\right)^{2/\alpha}{\pi}_{1}^{2}
,\ldots ,\left(\frac{n}{\log(n)}\right)^{2/\alpha}{\pi }_{m}^{2}\right)
^{T} \nonumber \\
{\longrightarrow} W_{1}^{2}+ \frac{m-1}{m} W^2_2+\ldots +\frac{1}{m}W_{m}^{2},
\newcounter{ATF}
\setcounter{ATF}{\value{equation}}
\end{eqnarray}
Finally, from eqn. ({\theATT}) and eqn. ({\theATF}),
\[
\left( 1,\frac{m-1}{m},\ldots,\frac{1}{m} \right)Y\rightarrow \sum^m_{i=1} \frac{m+1-i}{m} W^2_{i},
\]
and from ({\theAT}), we have the
\[
\hat D \rightarrow \sum^m_{i=1} \frac{m+1-i}{m} W^2_{i}. \ \ \ \ \ \ \Box
\]

\newpage
\bigskip
\begin{center}
APPENDIX B: MONTE-CARLO TEST PROCEDURE
\end{center}
\medskip

The Monte-Carlo test procedure for diagnostic checking of AR and ARMA models
with stable Paretian errors can be summarized below. Note that, to
check randomness of a time series, we skip Step 1 and in Step 4 we
simulate data from an IID sequence of $\{Z_{\hat \alpha}\}$ rather
than from the fitted model.

\begin{description}
\item[{\bf Step 1}] Fit an AR model to data using least-squares or the Burg algorithm or
for ARMA, an approximate Gaussian maximum likelihood algorithm is used.
Calculate residuals $\{\hat Z_t\}$ and the portmanteau
test of interest , say $\hat D_m$.
\item[{\bf Step 2}] Estimate $\alpha$ from residuals $\{\hat Z_t\}$ in Step 1.
The estimator given by McCulloch (1986) may be used.
\item[{\bf Step 3}] Select the number of Monte-Carlo simulations, $B$. Typically $100 \le B \le 1000$.
\item[{\bf Step 4}] Simulate the fitted model using the estimated AR or ARMA parameters in Step 1
and $\hat \alpha$ in Step 2.
Obtain $\hat D_m$ after estimating the parameters in the simulated series.
\item[{\bf Step 5}] Repeat Step 4 $B$ times counting the number of times $k$ that a value
of $\hat D_m$ greater than or equal to that in Step 1 has been
obtained.
\item[{\bf Step 6}] The $P$-value for the test is $(k+1)/(B+1)$.
\item[{\bf Step 7}] Reject the null hypothesis if the $P$-value is smaller than a predetermined significance level.
\end{description}

\newpage
\bigskip
\begin{center}
APPENDIX C: THE GENERALIZATION OF LINEAR EXPANSION OF RESIDUAL
AUTOCORRELATION
\end{center}
\medskip

\medskip
{\noindent \it C.1 Introduction \hfill}
\medskip

Residual autocorrelations are an important tool for diagnostic checking of
autoregressive and moving average ({\ARMA}) models. Their
asymptotic distributions from univariate {\ARMA} models were first
derived by Box and Pierce (1970). McLeod (1978) refined the
derivation and extended it to the multiplicative seasonal{\ARMA}
models. Their results were established under the assumption that
error sequences have finite variance and the parameters are
estimated using least squares, or equivalently, using maximum
likelihood estimation (MLE) for Gaussian {\ARMA} processes. Their
result may not be valid if the parameters of interest are estimated
using other estimation methods or linear processes with infinite
variance. This section demonstrates how the linear expansion of
residual autocorrelations in Box and Pierce (1970) also holds for
other estimation methods and for {\AR} models with stable Paretian
errors. The expansion may be used to derive the limiting
distribution of residual autocorrelations.

\medskip
{\noindent \it C.2 The Autoregressive Process \hfill}
\medskip

Consider an {\AR}$(p)$ process as follows:
 \begin{equation}
 \phi(B)y_t=a_t,
 \newcounter{ARP}
\setcounter{ARP}{\value{equation}}
 \end{equation}
where $B$ denotes the backward operator, $\phi(B)=1-\phi_1
B-\cdots-\phi_p B^p$, and $\{a_t\}$ is a sequence of independent and
identical random variables with mean zero and finite variance
$\sigma_a^2$. For given values $\dot{\Phi}=\left(\dot
\phi_1,\cdots,\dot \phi_p \right)^T $ of parameters, we can define
\begin{equation}
\dot a_t=a_t(\dot \Phi)=y_t-\dot \phi_1 y_{t-1}-\cdots-\dot \phi_p
y_{t-p}=\dot \Phi(B) y_t
\newcounter{innodot}
\setcounter{innodot}{\value{equation}}
\end{equation}
and the corresponding autocorrelation function at lag $k$ as
\begin{equation}
\dot r_k=r_k(\dot \Phi)=\frac{\sum \dot a_t \dot a_{t-k}}{\sum \dot
a^2_t}.
\end{equation}

\medskip
{\noindent \it C.3 Linear Expansion of Residual Autocorrelation
Function about Error Autocorrelation Functions \hfill}
\medskip

Consider approximating the residual autocorrelation $\hat r_k$
by a first order Taylor expansion about $\hat \Phi=\Phi$.
Let $\dot c_k$ and $\dot {r_k}$ denote $\sum \dot a_t \dot a_{t-k}$ and $\dot c_k/\dot c_0$
respectively, where $k \in$ integer.
Consider the estimators of $\Phi$ satisfying
\begin{equation}
\hat \phi_j=\phi_j+O_p\left(1/\sqrt{n}\right),\  \forall \ j.
\newcounter{phiorder}
\setcounter{phiorder}{\value{equation}}
\end{equation}
We have
\begin{equation}
\hat r_k=r_k+\sum^p_{j=1}\left( \phi_j-\hat \phi_j
\right)\hat \delta_{jk}+O_p\left( 1/n \right),
\newcounter{sacfexpand}
\setcounter{sacfexpand}{\value{equation}}
\end{equation}
where
\begin{eqnarray}
\hat \delta_{jk}
&=& -\frac{\partial \dot r_k}{\partial \dot \phi_j}|_{\dot\Phi=\hat \Phi} \nonumber \\
& = &
 -\frac{\partial}{\partial
\dot \phi_j}\left(\frac{\dot c_k}{\dot c_0}\right)|_{\dot\Phi=\hat \Phi} \nonumber \\
&=& \hat \delta^{(1)}_{ij}+\hat \delta^{(2)}_{ij},
\newcounter{rexpand}
\setcounter{rexpand}{\value{equation}}
\end{eqnarray}
\[\hat \delta^{(1)}_{ij}= -\dot c_k \frac{\partial }{\partial \dot \phi_j}
\left(\frac{1}{\dot c_0}\right)|_{\dot\Phi=\hat \Phi}
\]
 and
\[
\hat \delta^{(2)}_{ij}=
-\frac{1}{\dot c_0}
\frac{\partial \dot c_k}{\partial \dot \phi_j}|_{\dot\Phi=\hat \Phi}. \]
For {\rm LS} estimates, we have that
\begin{equation}
\frac{\partial}{\partial \dot \phi_j}\left[ \sum \dot a_t^2 \right]|_{\dot \Phi=\hat \Phi}
=\frac{\partial c_0}{\partial \dot \phi_j}|_{\dot \Phi=\hat \Phi}=0
\newcounter{lsresult}
\setcounter{lsresult}{\value{equation}}
\end{equation}
so it is straightforward that $\hat \delta_{ij}^{(1)}=0$.
Using this result, Box and Pierce (1970) showed that
$\hat \delta_{jk}=\psi_{k-j}$
to order $O_p\left(n^{-1/2}\right)$, where
$\psi_j$'s are the impulse response coefficients of the {\MA}$(\infty)$ representation
of eqn. ({\theARP}).
For other estimation methods, however,
$\hat \delta^{(1)}_{ij}$ may not be zero since eqn.
({\thelsresult}) does not hold.
To obtain a general result for $\hat \delta_{ij}$, therefore,
we will calculate $\hat \delta^{(1)}_{ij}$ explicitly.

Note that $\hat \delta^{(1)}_{ij}$ can be written as follows:
\begin{equation}
\dot c_k \cdot \left[ \sum \dot a_t^2\right]^{-2} \frac{\partial
\dot c_0}{\partial \dot \phi_j}|_{\dot \Phi=\hat \Phi}.
\newcounter{extraterm}
\setcounter{extraterm}{\value{equation}}
\end{equation}
By eqn. (2.15) of Box and Pierce (1970) and letting $k=0$, eqn. ({\theextraterm})
can be expressed as follows:
\begin{eqnarray}
& & \frac{\sum y_t^2}{\sum \hat{a_t^2}}
 \cdot \sum_{i=0}^p \hat
\phi_i
\left[r^{(y)}_{-i+j}+r^{(y)}_{i-j}\right]\cdot \frac{\hat c_k}{\hat c_0} \nonumber\\
 &=& \frac{\sum_{i=0}^p \hat
\phi_i
\left[r^{(y)}_{-i+j}+r^{(y)}_{i-j}\right]}{\sum_{i=0}^{p}\sum_{j=0}^{p}\hat
\phi_i \hat \phi_j r^{(y)}_{i-j}} \cdot \hat{r}_k,
\newcounter{secondterm}
\setcounter{secondterm}{\value{equation}}
\end{eqnarray}
where
\[
r_\nu^{(y)}=\frac{\sum y_t y_{t-\nu}}{\sum y_t^2}.
\]
Let $\hat \zeta_j$ denote
\[\left({\sum_{i=0}^p \hat \phi_i
\left[r^{(y)}_{-i+j}+r^{(y)}_{i-j}\right]}\right)/\left({\sum_{i=0}^{p}\sum_{j=0}^{p}\hat
\phi_i \hat \phi_j r^{(y)}_{i-j}}\right),
\]
and approximate $\hat \zeta_j$ by replacing $\hat \phi$'s and $r^{(y)}$'s
with $\phi$'s and $\rho$'s,
the theoretical parameters and the autocorrelations of the autoregressive process $\{y_t\}$.
By the Barteltt's formula,
\[
r_k^{(y)}=\rho_k+O_p\left(1/\sqrt n \right)
\]
as well as eqn. ({\thephiorder}) and ({\thesecondterm}), we have
\begin{equation}
\hat \zeta_j=\zeta_j
+O_p\left( 1/\sqrt{n} \right).
\end{equation}
Then by making use of the recursive relation
which is satisfied by the autocorrelations of an
autoregressive process,
eqn. (2.19) of Box and Pierce (1970), or
\begin{equation}
\rho_{\nu}-\phi_1 \rho_{\nu-1}-\cdots-\phi_p
\rho_{\nu-p}=\phi(B)\rho_{\nu}=0, \ \ \ \ \nu \geq 1,
\end{equation}
$\zeta_j$ can be simplified to yield
\begin{equation}
\zeta_j=\frac{\sum_{i=0}^p \phi_i \rho_{-j+i}}{\sum_{i=0}^{p}\phi_i
\rho_i}.
\newcounter{zetaform}
\setcounter{zetaform}{\value{equation}}
\end{equation}
Note that eqn. ({\thezetaform}) has the same form of eqn. (2.20) of
Box and Pierce (1970). Specifically, it can be seen as $\delta_{-j}$.
Moreover, Box and Pierce indicated that $\delta_{\nu}=0,\ \nu<0$ so $\zeta_j=0$.
Plugging this
result into eqn. ({\thesecondterm}), we have
$\hat \delta^{(1)}_{ij}=0 $. Consequently, eqn. (2.20) of Box and Pierce (1970)
for the linear expansion of residual autocorrelations
still holds for other estimators with order
$\hat \phi_i-\phi=O_p(1/\sqrt{n})$.

\medskip

{\noindent {\bf Remark 6 }: Many estimators of $\phi_{(p)}$ for an
{\AR} model with Paretian stable errors have order
$O_p([n/\log(n)]^{-1/\alpha})$, such as Whittle's, Yule-Walker and
{\rm LS} estimtors. Using the result that
$\textbf{r}_{(p)}=\rho_{(p)}+O_p([n/\log(n)]^{-1/\alpha})$, and
following the proofs in this section as well as in Box and Pierce
(1970), we may obtain the linear expansion of residual
autocorrelation functions for {\AR} models with stable Paretian
errors as in eqn. ({\thechthreeexpansion})}

\medskip
{\noindent \it C.4 The Equality of Residuals in {\AR} and {\ARIMA}
Models \hfill}
\medskip

The result in $\S C$.3 may be extended to {\ARIMA} models using
technique in $\S$5.1 of Box and Pierce (1970). If two time series
(a) an {\ARMA}($p,q$) process
\begin{equation}
\phi(B)w_t=\theta(B)a_t,
\newcounter{arma}
\setcounter{arma}{\value{equation}}
\end{equation}
and (b) an autoregressive series
\begin{equation}
\pi(B) x_t=\left(1-\pi_1 B- \cdots-\pi_{p+q}B^{p+q}\right)x_t=a_t,
\end{equation}
are both generated from the same set of errors $\{a_t\}$,
where
\[
\phi(B)=1-\phi B-\phi B^2-\cdots-\phi B^p,
\]
and
\[
\theta(B)=1-\theta B-\theta B^2-\cdots-\theta B^q.
\]
If
\begin{equation}
\pi(B)=\phi(B)\theta(B),
\newcounter{pieq}
\setcounter{pieq}{\value{equation}}
\end{equation}
then when the models are fitted by least squares, their residuals, and
hence also their autocorrelations, will be very nearly the same.
In this section, we consider whether the equality of residuals between
{\AR} and {\ARIMA} models is still valid when
the parameters are estimated by other approaches.

As in eqn. ({\theinnodot}), define
\begin{equation}
\dot a^{AR}_t=a_t^{AR}(\dot \pi)=\dot \pi(B)x_t=-\sum_{j=0}^{p+q}\dot \pi_j x_{t-j},
\newcounter{ARdot}
\setcounter{ARdot}{\value{equation}}
\end{equation}
where $\dot \pi_0=-1$, and now also
\begin{equation}
\dot a^\star_t=a_t^\star(\dot \phi,\dot \theta)=\dot \phi(B) \dot \theta(B)^{-1} w_t
=\left[\sum^p_{i=0} \dot \phi_i B^i\right] \left[\sum_{j=0}^q \dot \theta_j B^j\right]^{-1}w_t,
\newcounter{ARMAdot}
\setcounter{ARMAdot}{\value{equation}}
\end{equation}
where $\dot \phi_0=\dot \theta_0=-1$.
Using eqn. (5.12) and eqn. (5.13) of Box and Pierce (1970),
we can approximate $a^{AR}_t$ and $a^\star_t$ as follows:
\begin{equation}
\dot \mathbf{a}^{AR}=\mathbf{a}+ \mathbf{X} \left(\pi-\dot \pi \right)
\newcounter{ARapp}
\setcounter{ARapp}{\value{equation}}
\end{equation}
and
\begin{equation}
\dot \mathbf{a}^\star=\mathbf{a}+\mathbf{X} \left(\beta-\dot \beta\right).
\newcounter{ARMAapp}
\setcounter{ARMAapp}{\value{equation}}
\end{equation}

Note that eqn. ({\theARapp}) and eqn. ({\theARMAapp}) can be seen as a linear regression model.
We can estimate regression coefficients, $\pi-\dot \pi$ and $\beta-\dot \beta$ using
any suitable method. Let $g(\mathbf{X},\dot \mathbf{a}^\bullet)$ denote the corresponding estimator.
Since both eqn. ({\theARapp}) and eqn. ({\theARMAapp}) have the same form,
their estimators should agree with each other.
For example,
least squares estimates are given by
\begin{equation}
\hat \pi -\dot \pi=g(\mathbf{X},\dot \mathbf{a}^{AR})=(\mathbf{X}^T \mathbf{X})^{-1} \mathbf{X}^T \dot \mathbf{a}^{AR}
\newcounter{piest}
\setcounter{piest}{\value{equation}}
\end{equation}
and
\begin{equation}
\hat \beta -\dot \beta=g(\mathbf{X},\dot \mathbf{a}^\star)=(\mathbf{X}^T \mathbf{X})^{-1} \mathbf{X}^T \dot \mathbf{a}^\star.
\newcounter{betaest}
\setcounter{betaest}{\value{equation}}
\end{equation}
Then by setting $\dot \mathbf{a}=\mathbf{a}$
and estimating the regression coefficients of eqn. ({\theARapp}) and eqn. ({\theARMAapp}), we have
\begin{equation}
\hat \pi-\pi=g(\mathbf{X},\mathbf{a})=\hat \beta-\beta.
\newcounter{equalg}
\setcounter{equalg}{\value{equation}}
\end{equation}
Finally, by setting $\dot \mathbf{a}^{AR}= \hat \mathbf{a}^{AR}$ and
$\dot \mathbf{a}^\star= \hat \mathbf{a}^\star$
in eqn. ({\theARapp}) and eqn. ({\theARMAapp}),
it follows from
eqn. ({\theequalg}) that to order $O_p\left(|\hat\beta-\beta|^2\right)$
\begin{equation}
\hat \mathbf{a}^{AR}=g(\mathbf{X},\mathbf{a})=\hat \mathbf{a}^\star,
\end{equation}
and thus (to the same order) $\hat r^{AR}=\hat r^\star$.

\newpage
\begin{center}
{\smcaps REFERENCES\hfill}
\end{center}
\parindent 0pt

\hind
{\smcaps Adler, R.J. Feldman, R.E. and Gallagher, C.\/} (1998),
``Analysing Stable Time Series,''
{\it A Practical Guide to Heavy Tails: Statistical Techniques and Applications\/},
Birkh$\ddot{a}$user, Boston.

\hind
{\smcaps Box, G.E.P. and Pierce, D.A. \/} (1970),
``Distribution of Residual Autocorrelation in Autoregressive-Integrated Moving Average Time Series Models,''
{\it Journal of American Statistical Association\/} 65, 1509-1526.

\hind
{\smcaps Brockwell, P.J. and Davis, R.A. \/} (1991),
{\it Time Series: Theory and Methods\/},
Springer, New York.

\hind
{\smcaps Davis, R.A. \/} (1996),
``Gauss-Newton and M-estimation for {\ARMA} processes,''
{\it Stochastic Processes and their Applications\/} 63, 75--95.

\hind
{\smcaps Davis, R.A. and Resnick, S. \/} (1986),
``Limit Theory for the Sample Covariance and Correlation Functions of Moving Averages,''
{\it The Annals of Statistics\/} 14, 533--558.



\hind {\smcaps McCulloch, J.H. \/} (1986), ``Simple Consistent
Estimators of Stable Distribution Parameters,'' {\it Communication
in Statistics--Computation and Simulation\/}, {15}, 1109--1136.

\hind
{\smcaps Mittnik, S., Rachev, S.T., Doganoglu, T. and Chenyao, D. \/} (1999),
``Maximum Likelihood Estimation of Stable Paretian Models,''
{\it Mathematical and Computer Modelling\/}, {29}, 275--293.

\hind
{\smcaps Monti, A.C. \/} (1994),
``A Proposal for Residual Autocorrelation Test in Linear Models,''
{\it Biometrika\/} 81, 776--780.

\hind
{\smcaps Pe\v{n}a, D. and Rodriguez, J. \/} (2002),
``A Powerful Portmanteau Test of Lack of Fit For Time Series,''
{\it Journal of American Statistical Association\/} 97, 601-610.

\hind
{\smcaps Runde, R. \/} (1997),
``The Asymptotic Null Distribution of the Box-Pierce Q-Statistic for Random Variable with Infinite Variance:
 An Application to German Stock Returns,''
{\it Journal of Econometrics\/} 78, 205-216.

\hind
{\smcaps  Samorodnitsky, G. and Taqqu, M. \/} (1994),
{\it Stable-Non-Gaussian Random Processes\/},
Chapman-Hall, New York.

\hind
{\smcaps Tsay, R.S.\/} (2002),
 {\it Analysis of Financial Time
Series\/}, New York: Wiley.

\newpage
\begin{table}
{\smcaps Table I.
Empirical sizes $(\%)$ of $\hat D$ and $Q_{\rm BP}$
for a $5\%$ significance test
based on the parametric bootstrap procedure.
The empirical size for each test was
calculated based on $N=10^4$ simulations.
Each Monte Carlo test also used $B=10^3$ simulations.
Series length $n=250$ and lags $m=5,10,15$ were investigated.
}
\begin{center}
\begin{tabular}{ccccccc}
\noalign{\smallskip}
\noalign{\hrule}
\noalign{\smallskip}
\noalign{\hrule}
\noalign{\smallskip}
\noalign{\smallskip}
$     $ & $\hat D(5) $ &$\hat D(10)  $  &$\hat D(15)$ & $Q_{\rm BP}(5)$ &$Q_{\rm BP}(10)$  &$Q_{\rm BP}(15)$     \\
$\alpha=1.9 $ & $5.30$    &$4.66$     &$4.78$ & $4.96$    &$4.71$     &$4.87$        \\
$\alpha=1.7$ & $5.18$    &$4.44$     &$4.44$ & $4.82$    &$4.43$     &$4.41$        \\
$\alpha=1.5 $ & $4.82$    &$4.99$     &$5.13$ & $5.07$    &$5.27$     &$5.30$        \\
$\alpha=1.3$ & $4.80$    &$5.03$     &$5.18$ & $5.04$    &$5.00$     &$5.27$        \\
$\alpha=1.1$ & $5.26$    &$5.33$     &$5.12$ & $5.33$    &$5.25$     &$5.15$        \\
\noalign{\smallskip}
\noalign{\hrule}
\end{tabular}
\end{center}
\end{table}
\strut


\strut \eject\vfill\clearpage \newpage
\begin{table}
{\smcaps Table II.
P-values for $Q_{\rm LB}$ statistic using Monte-Carlo test and
$\chi^2$-method for testing randomness of exchange-rate returns.
 }
\begin{center}
\begin{tabular}{lcc}
\noalign{\smallskip} \noalign{\hrule} \noalign{\smallskip}
\noalign{\hrule} \noalign{\smallskip}
\noalign{\smallskip}
$             $ & Monte-Carlo Test & $\chi^2(m)$ Test\\
$m=5$          &$0.500$     &$0.042$       \\
$m=10$         &$0.582$     &$0.228$      \\
$m=20$         &$0.828$     &$0.404$     \\
\noalign{\smallskip} \noalign{\hrule}
\end{tabular}
\end{center}
\end{table}

\newpage
\begin{table}
{\smcaps Table III.
Empirical sizes $(\%)$ of $\hat D$ and $Q_{\rm BP}$
for a $5\%$ significance test.
$\hat D$ and $Q_{\rm BP}$ tests for checking model adequacy of
{\AR}$(1)$ models
fitted by the Burg algorithm.
Both tests
were implemented by the parametric bootstrap procedure.
The empirical size for each test was
calculated based on $N=10^4$ simulations.
Each Monte Carlo test also used $B=10^3$ simulations.
Series length $n=100$ and lags $m=5,10,20$ were investigated.
}
\begin{center}
\begin{tabular}{rcccccc}
\noalign{\smallskip}
\noalign{\hrule}
\noalign{\smallskip}
\noalign{\hrule}
\noalign{\smallskip}
\noalign{\smallskip}
$\phi_1     $ & $\hat D(5) $ &$\hat D(10)  $  &$\hat D(20)$ & $Q_{\rm BP}(5)$ &$Q_{\rm BP}(10)$  &$Q_{\rm BP}(20)$     \\
$0.9$  & $4.90$    &$4.75$     &$4.88$ & $4.60$    &$4.71$     &$4.96$        \\
$0.7$  & $4.97$    &$5.20$     &$5.16$ & $4.95$    &$4.94$     &$5.42$        \\
$0.5$  & $5.37$    &$5.32$     &$5.14$ & $5.55$    &$5.12$     &$5.16$        \\
$0.3$  & $5.11$    &$4.90$     &$4.82$ & $5.13$    &$4.80$     &$5.26$        \\
$0.1$  & $4.92$    &$5.01$     &$5.20$ & $5.14$    &$4.75$     &$4.86$        \\
$-0.1$ & $5.30$    &$5.45$     &$5.29$ & $5.25$    &$5.08$     &$4.90$        \\
$-0.3$ & $5.00$    &$5.20$     &$5.33$ & $4.79$    &$5.30$     &$5.45$        \\
$-0.5$ & $5.00$    &$4.93$     &$5.10$ & $5.00$    &$4.93$     &$5.26$        \\
$-0.7$ & $5.62$    &$5.73$     &$5.65$ & $5.20$    &$5.45$     &$5.41$        \\
$-0.9$ & $5.21$    &$5.02$     &$5.07$ & $5.01$    &$5.00$     &$5.30$        \\
\noalign{\smallskip}
\noalign{\hrule}
\end{tabular}
\end{center}
\end{table}
\strut
\vfill\eject\newpage

\begin{table}
{\smcaps Table IV.
Empirical powers $(\%)$ of $\hat D$ and $Q_{\rm BP}$
for a $5\%$ significance test.
$\hat D$ and $Q_{\rm BP}$ tests for checking model adequacy of
twelve {\ARMA}$(2,2)$ models in Table 3 of
Pe\v{n}a and Rodriguez (2002)
fitted by {\AR}$(1)$ using the Burg algorithm.
Both tests
were implemented based on the parametric bootstrap procedure.
The empirical power for each test was
calculated based on $N=10^4$ simulations.
Each Monte Carlo test also used $B=10^3$ simulations.
Series length $n=100$ and lags $m=5,10,20$ were investigated.
}
\begin{center}
\begin{tabular}{ccccccc}
\noalign{\smallskip}
\noalign{\hrule}
\noalign{\smallskip}
\noalign{\hrule}
\noalign{\smallskip}
\noalign{\smallskip}
 Model      & $\hat D(5) $ &$\hat D(10)  $  &$\hat D(20)$ & $Q_{\rm BP}(5)$ &$Q_{\rm BP}(10)$  &$Q_{\rm BP}(20)$     \\
1  & $53.32$    &$38.31$     &$32.77$ & $29.59$    &$21.76$     &$19.25$        \\
2  & $99.01$    &$98.56$     &$98.01$ & $94.53$    &$70.46$     &$59.61$        \\
3  & $23.37$    &$20.10$     &$17.61$ & $21.62$    &$16.71$     &$15.17$        \\
4  & $77.13$    &$59.38$     &$48.12$ & $60.82$    &$40.29$     &$35.15$        \\
5  & $93.22$    &$87.62$     &$79.84$ & $84.66$    &$66.68$     &$58.46$        \\
6  & $13.74$    &$11.17$     &$10.05$ & $10.68$    &$9.13$      &$8.61$        \\
7  & $26.51$    &$26.25$     &$24.92$ & $17.56$    &$13.80$     &$13.05$        \\
8  & $33.92$    &$26.68$     &$23.57$ & $27.36$    &$20.60$     &$19.25$        \\
9  & $99.44$    &$99.27$     &$99.16$ & $98.71$    &$93.17$     &$78.88$        \\
10 & $76.71$    &$58.06$     &$48.50$ & $40.62$    &$28.39$     &$25.94$        \\
11 & $99.01$    &$98.46$     &$97.87$ & $94.02$    &$67.04$     &$57.11$        \\
12 & $99.89$    &$99.87$     &$99.48$ & $99.86$    &$99.63$     &$99.48$        \\
\noalign{\smallskip}
\noalign{\hrule}
\end{tabular}
\end{center}

\end{table}

\strut \eject\vfill\clearpage \newpage
\begin{table}
{\smcaps Table V. An illustrated example using the monthly simple
return of {\rm CRSP} value-weighted index data from Tsay (2002). The
data were fitted by an AR$(3)$ model and an AR$(5)$ model. The
entries in the first two columns are the P-values of $\hat D$ and
$Q_{\rm BP}$ in \S4 based on the Monte-Carlo test; those in the
third column are the P-value of the portmanteau test of Box and
Pierce (1970) assuming a normal distribution, denoted by $Q^{\rm
N}_{\rm BP}$. }
\begin{center}
\begin{tabular}{cccc}
\noalign{\smallskip} \noalign{\hrule} \noalign{\smallskip}
\noalign{\hrule} \noalign{\smallskip}
&\multicolumn{3}{c}{{\AR}$(3)$} \\
\noalign{\smallskip}
$             $ & $\hat D$ &$Q_{\rm BP}$  &$Q_{\rm BP}^{\rm N}$\\
$m=5$        & $0.050$    &$0.026$     &$0.197$       \\
$m=10$        & $0.030$    &$0.021$     &$0.107$      \\
$m=20$        & $0.019$    &$0.012$     &$0.247$     \\
\noalign{\smallskip} \noalign{\hrule} \noalign{\smallskip}
\noalign{\smallskip}
&\multicolumn{3}{c}{{\AR}$(5)$} \\
\noalign{\smallskip}
$             $ & $\hat D$ &$Q_{\rm BP}$  &$Q_{\rm BP}^{\rm N}$\\
$m=5$         & $0.064$     &$0.055$     &$0.998$       \\
$m=10$        & $0.052$    &$0.045$     &$0.345$       \\
$m=20$        & $0.024$    &$0.024$     &$0.438$       \\
\noalign{\smallskip} \noalign{\hrule}
\end{tabular}
\end{center}
\end{table}

\newpage
\begin{figure}[h]
\centerline{\epsfig{figure=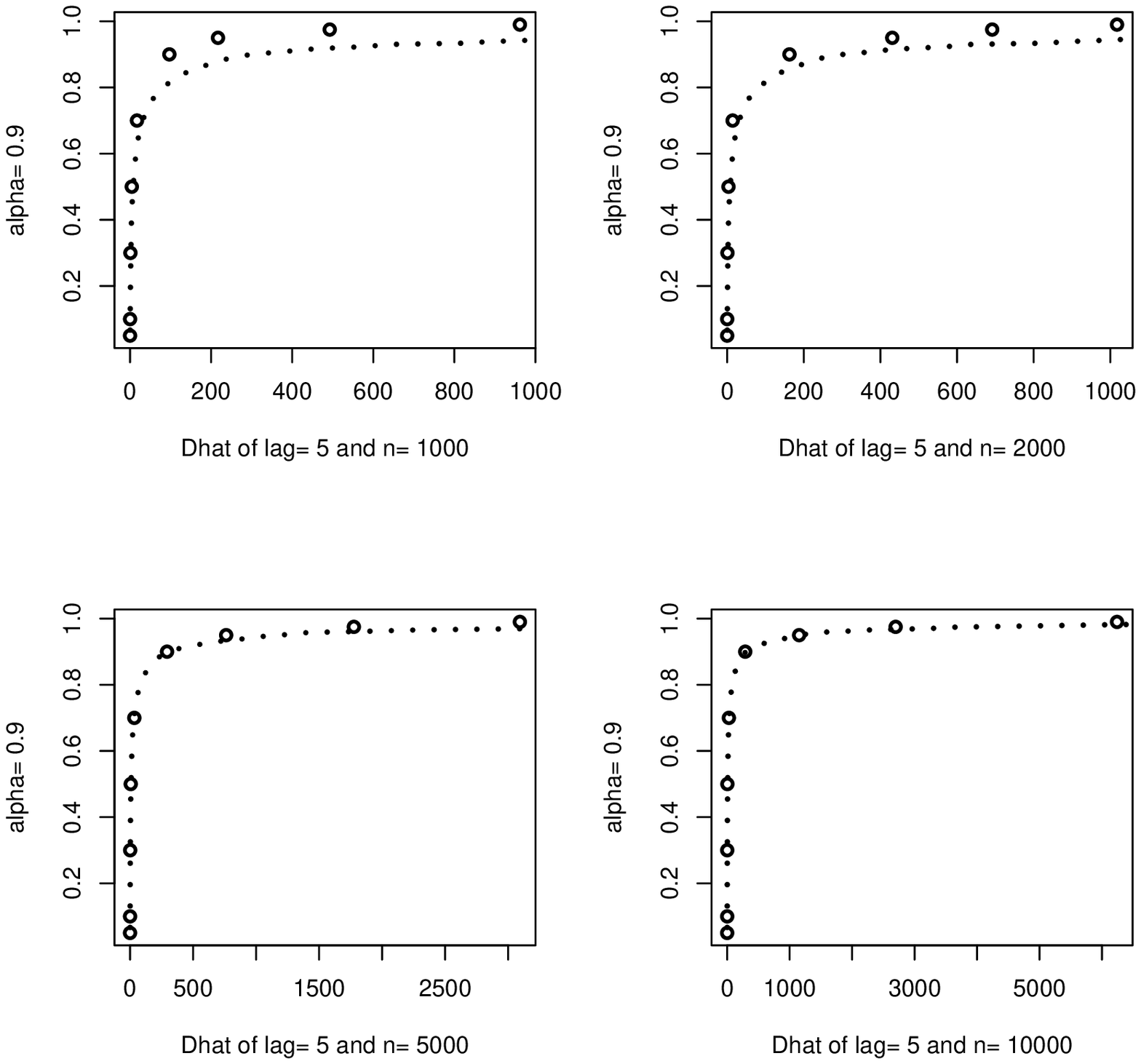, width=6.5 in,height=5.5 in}}
\caption{
The slow convergence of the $\hat D$ test to its asymptotic distribution.
Random sequences of series length $n=10^3,2000,5000,10^4$ were simulated from $S_{1.5}(1,0,0)$.
$250$ simulations were used to
retrieve empirical percentiles of the $\hat D$ test with $m=5$.
The $5$, $10$, $30$, $50$, $70$, $90$, $95$, $97.5$, $99$ $(\%)$ empirical quantiles
were plotted as black circles and
the corresponding asymptotic distribution was also plotted as the dot line.}
\end{figure}

\newpage
\begin{figure}[h]
\centerline{\epsfig{figure=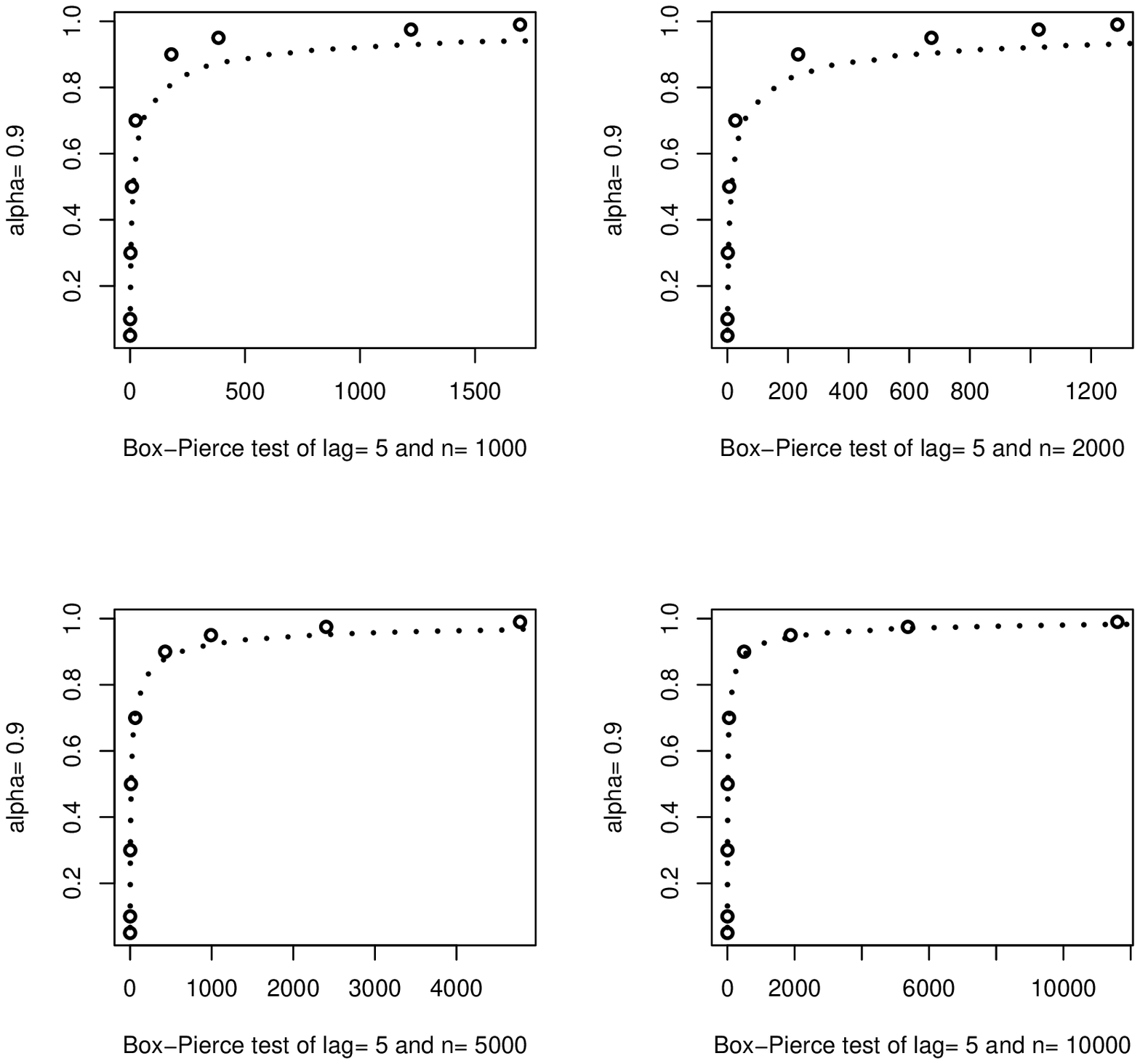, width=6.5 in,height=5.5 in}}
\caption{
The slow convergence of the $Q_{\rm BP}$ test to its asymptotic distribution.
Random sequences of series length $n=10^3,2000,5000,10^4$ were simulated from $S_{1.5}(1,0,0)$.
$250$ simulations were used to
retrieve empirical percentiles of the $Q_{\rm BP}$ test with $m=5$.
The $5$, $10$, $30$, $50$, $70$, $90$, $95$, $97.5$, $99$ $(\%)$ empirical quantiles
were plotted as circles and
the corresponding asymptotic distribution was also plotted as the dot line.}
\end{figure}

\newpage
\begin{figure}[h]
\centerline{\epsfig{figure=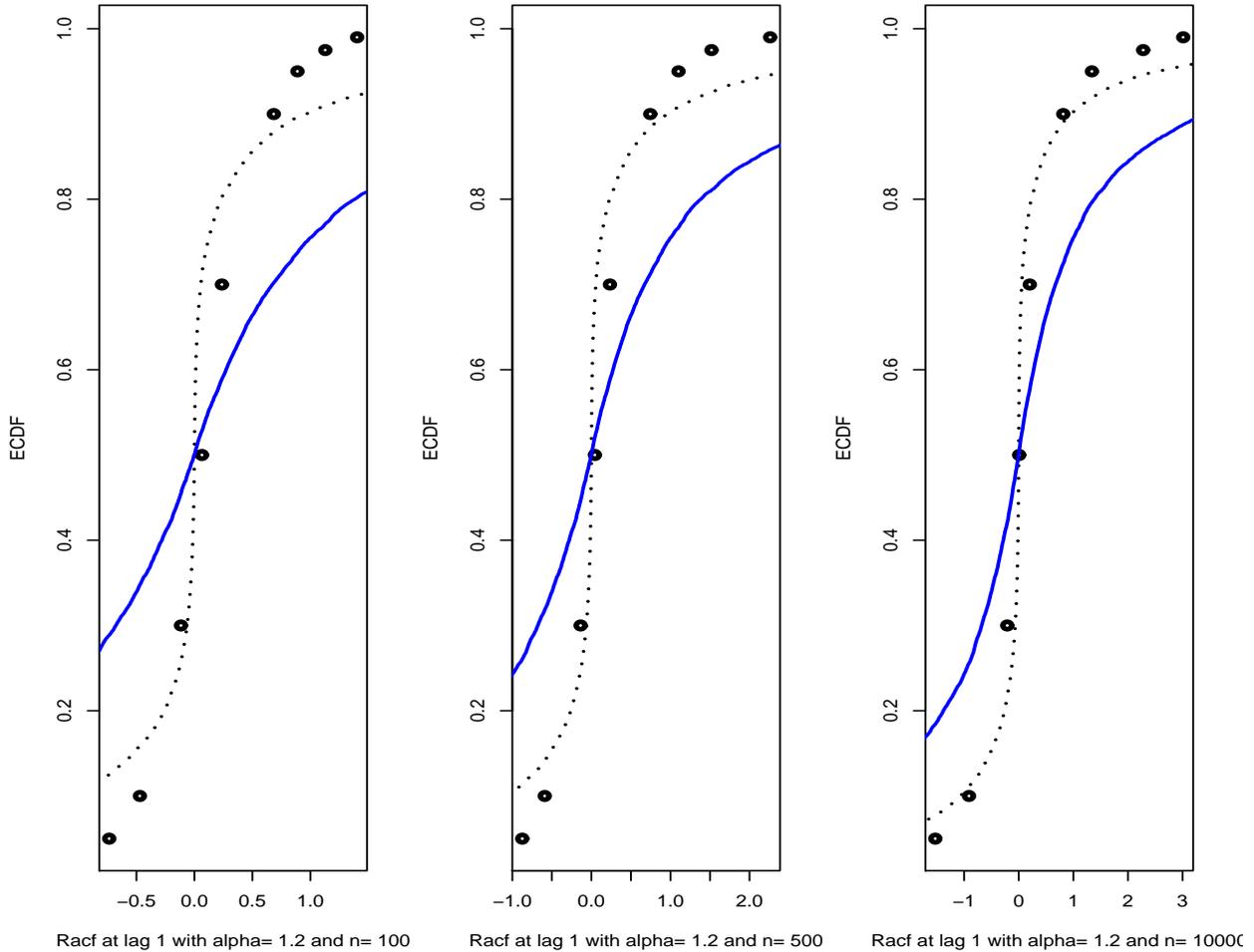, width=6.5 in,height=5.5 in}}
\caption{ The slow convergence of residual autocorrelation to its
asymptotic distribution.{\AR}$(1)$ process, $X_t=0.5X_{t-1}+Z_t$, of
series length $n=100$, $500$, $10^4$ were simulated respectively,
where $\{Z_t\}$ is distributed as $Z_{1.2}(1,0,0)$. The number of
simulation ${\rm ~NSIM}=10^4$ were used.{\AR}$(1)$ models were then
fitted to simulated data and residual autocorrelation at lag one was
calculated. The $5$, $10$, $30$, $50$, $70$, $90$, $95$, $97.5$,
$99$ $(\%)$ empirical quantiles of residual autocorrelation at lag
one were plotted as circles. The corresponding asymptotic
distribution was plotted as the dot line. The asymptotic
distribution of sample autocorrelation was plotted as the real line.
}
\end{figure}

\end{document}